\documentclass[11pt]{article}

\usepackage{amsmath,amsthm,enumitem,fontspec,unicode-math,geometry,hyperref}
\providecommand{\restriction}{\mathord{↾}}
\newcommand{\pre}{\mathrel{⊏}}
\hypersetup{hidelinks,pdftitle={A topological characterization of universally Baire sets in ZF},pdfauthor={Zihao Guo}}
\setlist{nosep,leftmargin=2em}
\allowdisplaybreaks

\newtheorem{theorem}{Theorem}[section]
\newtheorem{lemma}[theorem]{Lemma}

\newtheorem{corollary}[theorem]{Corollary}
\theoremstyle{definition}

\newtheorem{remark}[theorem]{Remark}

\newcommand{\R}{\omega^\omega}
\newcommand{\RO}{\operatorname{RO}}
\newcommand{\BA}{\operatorname{BA}}
\newcommand{\Pow}{\mathcal P}
\newcommand{\cl}{\operatorname{cl}}
\newcommand{\intt}{\operatorname{int}}

\title{A topological characterization of universally Baire sets in ZF}
\author{Zihao Guo}
\date{}

\begin{document}
\begin{center}
{\Large A topological characterization of\\[2pt] universally Baire sets in ZF\par}
\vspace{7pt}Zihao Guo
\end{center}

\begin{abstract}
In ZF we prove that a set of reals whose continuous preimages have the Baire property in all spaces with a regular open base admits absolutely complementing ordinal trees for each set forcing. This gives the converse of a lemma of Larson, Sargsyan, and Wilson. For each forcing we construct one completely ultrametrizable test space, independently of the set of reals. A Baire decomposition on this space determines families of finite extensions in a complete Boolean algebra. Whenever the supremum of such a family leaves a nonzero remainder, that remainder lies outside the current subalgebra. A Hartogs recursion then produces a well-orderable subalgebra in which all the required suprema are exact. Its finite histories supply the labels of both trees.

\end{abstract}


\section{Introduction}

We call the elements of Baire space $\omega^\omega$ reals. An open set is regular open if it is the interior of its closure. For a set of reals $A$, consider the following condition.
\begin{quote}
$(\ast)$ For every space $X$ with a regular open base and every continuous $f:X\to\omega^\omega$, the preimage $f^{-1}[A]$ has the Baire property in $X$.
\end{quote}

A set has the Baire property if its symmetric difference with an open set is contained in a countable union of closed nowhere dense sets. In ZF, allowing arbitrary nowhere dense sets gives the same definition, by replacing each such set with its closure.

For an ordinal $\theta$ and a tree $T\subseteq(\omega\times\theta)^{<\omega}$, put
\[
p[T]=\{x\in\R:\exists z\in\theta^\omega\ \forall n<\omega\ ((x\restriction n,z\restriction n)\in T)\}.
\]
A pair $(T^+,T^-)$ of trees on a common ordinal alphabet is $\mathbb P$-absolutely complementing if
\[
1_{\mathbb P}\Vdash p[\check T^+]=\omega^\omega\setminus p[\check T^-].
\]
The projections in this assertion are computed in the forcing extension, using all branches in that extension.

\begin{theorem}\label{thm:main}
In ZF the following are equivalent for $A\subseteq\R$.
\begin{enumerate}[label=(\roman*)]
\item $A$ satisfies $(\ast)$.
\item For every set forcing $\mathbb P$ there are an ordinal $\theta$ and a $\mathbb P$-absolutely complementing pair $T^+,T^-$ on $\omega\times\theta$ such that
\[
p[T^+]^V=A,\qquad p[T^-]^V=(\R)^V\setminus A.
\]
\end{enumerate}
\end{theorem}

Feng, Magidor, and Woodin established this equivalence in ZFC \cite[Corollary 2.1(3)]{FMW}. In ZF, Larson, Sargsyan, and Wilson prove the implication from trees to $(\ast)$. Immediately before Lemma 2.4 they state that they do not know whether the converse is provable in ZF \cite{LSW}. We prove this converse. In their ZFC comparison of test spaces, Elekes and P\'alfy show that arbitrary spaces, spaces with a regular open base, and completely metrizable spaces give equivalent tests \cite[Claim~4.4, parts (1), (3), and (6)]{EP}. Section~5 derives a complete ultrametric test in ZF from the particular spaces used here. A related characterization of Ikegami and Viale concerns a fixed complete Boolean algebra $B$ and subsets of $2^\kappa$. Under $\mathrm{FA}_\kappa(B\restriction b)$ for each nonzero $b$, their Theorem 14 allows trees on $2\times Y$ for a set $Y$ \cite{IV}. Here the auxiliary alphabet must be well-orderable even when the forcing and its set of real names need not be.

\paragraph{Why measure a supremum?}
A Baire decomposition can assign a color correctly on a comeager set and still give the wrong answer at the limit of a particular chain. For example, the singleton $\{0^\omega\}$ is nowhere dense, yet every Cohen condition $0^n$ is compatible with that real. We therefore require a branch to enter a prescribed sequence of dense open sets. This ensures correctness once a branch exists. It leaves a second question: after imposing the next requirement, is there still an extension reading the prescribed real?

The distinction already occurs in Cantor space. The cylinders $[0^n1]$ have supremum $1$ in the algebra of clopen sets, since their union is dense. Their supremum in the full power-set algebra misses $0^\omega$. Measuring the join in the larger algebra detects precisely the lost point. In our construction the ambient algebra $D$ has a coordinate for every old real and another for every Boolean code for a real. An equality of joins in $D$ can thus be tested at any of these coordinates.

We obtain those equalities before coding histories by ordinals. First we apply $(\ast)$ to the topological sum of the history spaces of all candidate subalgebras. One Baire decomposition gives dense requirements on every summand. Whenever the required join falls short in $D$, its remainder is an element outside the current subalgebra. Adjoining the remainder for the first failed requirement cannot continue up to the Hartogs ordinal of $D$. At the stopping stage the histories are well-orderable and their extensions cover each prescribed coordinate. Section~4 uses this last assertion twice: directly for an old real, and through genericity for a real in an extension.

\section{Boolean-valued coding and the test space}

Fix a nonempty set forcing $\mathbb P$. If necessary, adjoin a largest condition; this does not change its generic extensions. We use the convention that smaller conditions are stronger. Let
\[
\mathbb B=\RO(\mathbb P)
\]
be the complete Boolean algebra of regular-open subsets of $\mathbb P$, where the topology consists of downward-closed sets. Write \(\mathbb B^+=\mathbb B\setminus\{0\}\) for the positive part of \(\mathbb B\). For $p\in\mathbb P$, write $\downarrow p=\{q\in\mathbb P:q\leq p\}$ and let $\rho(p)=\intt\cl(\downarrow p)$. For any family $\mathcal U\subseteq\RO(\mathbb P)$ its join is $\intt\cl(\bigcup\mathcal U)$, the least regular open set containing its union. Boolean complementation sends $U$ to $\mathbb P\setminus\cl U$. These operations exist by Separation on the power set of $\mathbb P$. The map $p\mapsto\rho(p)$ has dense range in $\mathbb B^+$: if $U$ is a nonempty regular-open set and $p\in U$, then $\rho(p)\subseteq U$. Membership in the regularization has the following concrete description:
\[
 q\in\rho(p)\quad\Longleftrightarrow\quad
 (\forall r\le q)(\exists s\le r)\ s\le p.
\]
Indeed, $\downarrow q$ is the smallest neighborhood of $q$, and $r$ belongs to the closure of $\downarrow p$ exactly when $r$ and $p$ have a common extension. Thus $\rho(p)\cap\rho(p')$ is nonempty exactly when $p$ and $p'$ are compatible. For the forward implication, take $q$ in the intersection, refine $q$ below $p$, and then refine once more below $p'$. Conversely a common extension belongs to both regularizations. The map therefore gives the usual dense embedding of the separative quotient, so forcing with $\mathbb P$ is equivalent to forcing with $\mathbb B^+$; we work with the latter. A $\mathbb B^+$-generic filter is identified with its induced ultrafilter on $\mathbb B$.

\subsection{Boolean matrices for reals}

Let $I$ be the set of all matrices $i=(i_{n,m})_{n,m<\omega}\in\mathbb B^{\omega\times\omega}$ such that, for each $n<\omega$,
\[
i_{n,m}\wedge i_{n,k}=0\quad(m\ne k),\qquad \bigvee_{m<\omega}i_{n,m}=1.
\]
Zero entries are allowed. If $G$ is $\mathbb B^+$-generic, define
\[
i_G(n)=m\quad\Longleftrightarrow\quad i_{n,m}\in G.
\]
For each fixed $n$, the nonzero entries in row $n$ form a maximal antichain: their join is $1$. Hence $G$ meets this antichain, and pairwise disjointness ensures that it meets at most one member. Thus $i_G\in\R$ is well-defined.

\begin{lemma}\label{lem:all-reals}
Every real in $V[G]$ is equal to $i_G$ for some $i\in I^V$.
\end{lemma}
\begin{proof}
Let $x\in(\R)^{V[G]}$, and choose a ground-model name $\sigma$ whose interpretation is the graph of $x$. By the truth lemma, there is $b\in G$ such that
\[
b\Vdash_{\mathbb B}\text{``$\sigma$ is a function from $\check\omega$ to $\check\omega$.''}
\]
For $n,m<\omega$, put
\[
c_{n,m}=b\wedge\bigl\|\langle\check n,\check m\rangle\in\sigma\bigr\|_{\mathbb B}.
\]
For fixed $n$, the elements $c_{n,m}$ are pairwise disjoint, since $b$ forces uniqueness of the value of $\sigma(n)$. Their join is $b$: below $b$, conditions deciding the natural-number value of $\sigma(n)$ are dense. Define
\[
i_{n,0}=c_{n,0}\vee\neg b,\qquad i_{n,m}=c_{n,m}\quad(m>0).
\]
Each row is now a partition of unity, so $i\in I^V$. Since $b\in G$, the added piece $\neg b$ is not in $G$, and the evaluation of the matrix agrees with $\sigma_G$ at every coordinate. Hence $i_G=x$.
\end{proof}

The given name need only be forced to be a real below $b$. Assigning the fixed value $0$ on $\neg b$ makes the code total. This argument uses neither the maximum principle nor a simultaneous choice of names.

\subsection{The ambient Boolean algebra}

Let $R$ denote the ground-model set $(\R)^V$, and form the set-sized complete Boolean algebra
\[
D=\Pow(R)\times\mathbb B^I,
\]
with coordinatewise operations. Explicitly, for $H\subseteq D$ its join is
\[
 \left(\bigcup_{(E,f)\in H}E,\quad i\longmapsto\bigvee_{(E,f)\in H} f(i)\right).
\]
For each $i$ the displayed value is unique, so Replacement gives the function on all of $I$. The least-upper-bound property follows in each coordinate. For $n,m<\omega$, define
\[
d_{n,m}=\bigl(\{x\in R:x(n)=m\},\ (i_{n,m})_{i\in I}\bigr).
\]
For each fixed $n$, the family $(d_{n,m})_{m<\omega}$ is a partition of unity in $D$.

We use two kinds of complete Boolean homomorphisms out of $D$. For $x\in R$, let
\[
\pi_x:D\to 2,\qquad \pi_x(E,f)=1\ \Longleftrightarrow\ x\in E.
\]
For $i\in I$, let
\[
\pi_i:D\to\mathbb B,\qquad \pi_i(E,f)=f(i).
\]
Both maps preserve arbitrary joins and complements. The first is evaluation at a ground-model real in the first factor; the second is projection to one Boolean-name coordinate in the second factor. These two maps will later provide, respectively, ground-model coverage and coverage in a generic extension.

\subsection{Histories and their topological spaces}

Let $\mathcal C$ be the set of all Boolean subalgebras $C\subseteq D$ containing every $d_{n,m}$. Only closure under the finite Boolean operations is required. For $C\in\mathcal C$, a \emph{history} is a finite sequence
\[
t=(a_0,\ldots,a_{k-1})
\]
of nonzero elements of $C$ such that
\[
a_0\geq a_1\geq\cdots\geq a_{k-1}
\]
and, for each $j<k$, there is an $m_j<\omega$ with $a_j\leq d_{j,m_j}$. Equality among successive $a_j$ is allowed. The value $m_j$ is unique by disjointness. Write $a_t=a_{k-1}$ and $u_t=(m_0,\ldots,m_{k-1})$. For the empty history put $a_\varnothing=1$ and $u_\varnothing=\varnothing$. Let $L_C$ be the set of histories.

\begin{lemma}\label{lem:extend}
If $t\in L_C$ and $0<b\leq a_t$ with $b\in C$, then there is a one-step extension $v$ of $t$ such that $0<a_v\leq b$. In particular, every history lies on an infinite history.
\end{lemma}
\begin{proof}
Let $k=|t|$. Since $(d_{k,m})_{m<\omega}$ is a partition of unity and meet distributes over arbitrary joins in $D$,
\[
b=b\wedge 1=b\wedge\bigvee_m d_{k,m}=\bigvee_m(b\wedge d_{k,m}).
\]
As $b>0$, some $m$ satisfies $b\wedge d_{k,m}>0$. Choose the least such $m$ and append $b\wedge d_{k,m}$ to $t$. This element belongs to $C$ by finite closure, and it is below $b$. Repeating this uniquely specified rule with the current last element in place of $b$ gives an infinite history by recursion on $\omega$. No choice principle is needed.
\end{proof}

Let
\[
X_C=\{y\in C^\omega:\ y\restriction k\in L_C\text{ for every }k<\omega\}.
\]
Give $C$ the discrete topology and $X_C$ the subspace topology inherited from $C^\omega$. For $t\in L_C$, put
\[
[t]_C=\{y\in X_C:t\subseteq y\}.
\]
By Lemma~\ref{lem:extend}, every such cylinder is nonempty. It is clopen, since it specifies finitely many coordinates in a product of discrete spaces. These cylinders form a base: any finite-coordinate neighborhood of $y$ contains the cylinder determined by a sufficiently long initial segment of $y$. Thus $X_C$ has a clopen base.

For $y\in X_C$, define $e_C(y)\in\R$ by
\[
e_C(y)(j)=m\quad\Longleftrightarrow\quad y(j)\leq d_{j,m}.
\]
The history condition and disjointness make this well-defined. The first $k$ output coordinates are constant on each length-$k$ cylinder, so $e_C:X_C\to\R$ is continuous.

The later recursion determines its next algebra from the Baire data on the current algebra. Applying $(\ast)$ separately at each stage would give individual witnesses without specifying a sequence of witnesses for this recursion. To obtain all the data at once, take the topological sum
\[
X=\coprod_{C\in\mathcal C}X_C,
\qquad e(C,y)=e_C(y).
\]
Each summand is clopen. The tagged cylinders from all summands together form a clopen base for $X$. In particular, $X$ has a regular-open base and $e:X\to\R$ is continuous. The space $X$ and the map $e$ depend on $\mathbb P$ but not on $A$.

\section{Ambient suprema and the Hartogs refinement}

Assume now that $A\subseteq\R$ satisfies $(\ast)$. Applying $(\ast)$ to the single map $e:X\to\R$, choose an open set $U\subseteq X$ and closed nowhere dense sets $F_n\subseteq X$ such that
\[
e^{-1}[A]\mathbin{\triangle} U\subseteq\bigcup_{n<\omega}F_n.
\]
Put
\[
U_+=U,\qquad U_-=X\setminus\cl(U),\qquad O_n=X\setminus F_n,
\]
and let $A_+=A$ and $A_-=\R\setminus A$. The two color sets $U_+$ and $U_-$ are disjoint and open, and their union is
\[
U_+\cup U_-=X\setminus\partial U,
\]
a dense open set. Each $O_n$ is also dense open. On ground-model points we therefore have
\begin{equation}\label{eq:color-correct}
y\in U_\varepsilon\cap\bigcap_{n<\omega}O_n
\quad\Longrightarrow\quad e(y)\in A_\varepsilon.
\end{equation}
Indeed, a point meeting every $O_n$ is outside the symmetric difference. In $U_+$ its image therefore belongs to $A$. In $U_-$ the point lies outside $U$, so its image belongs to the complement of $A$.

In cylinder inclusions we identify $X_C$ with its tagged copy in $X$. Write $t\pre s$ when $t$ is a proper initial segment of $s$. For $t\in L_C$ and $n<\omega$, define
\[
J_C(t,n)=\{s\in L_C:t\pre s,\ (\exists\varepsilon\in\{+,-\})\ [s]_C\subseteq O_n\cap U_\varepsilon\},
\]
and take the join in the ambient algebra $D$:
\[
j_C(t,n)=\bigvee_D\{a_s:s\in J_C(t,n)\}.
\]
Every contributing $a_s$ is below $a_t$, so $j_C(t,n)\leq a_t$. Define the residual
\[
r_C(t,n)=a_t\wedge\neg j_C(t,n).
\]
The join is deliberately taken in $D$; it need not belong to $C$.

\begin{lemma}[Internal supremum and ambient remainder]\label{lem:residual}
Fix $C\in\mathcal C$, $t\in L_C$, and $n<\omega$, and set $H=\{a_s:s\in J_C(t,n)\}$. Then $\sup_C H$ exists and equals $a_t$. If $r_C(t,n)>0$, then $r_C(t,n)\notin C$.
\end{lemma}
\begin{proof}
We first show that every nonzero $b\in C$ below $a_t$ bounds some member of $H$ from above. By Lemma~\ref{lem:extend}, extend $t$ to $v$ with $0<a_v\le b$. The nonempty cylinder $[v]_C$ meets $O_n\cap(U_+\cup U_-)$, a dense open set on this summand. At a point of the intersection choose an initial segment $s$ extending $v$ whose cylinder lies in $O_n$ and in the color open set containing the point. Then $s\in J_C(t,n)$ and $0<a_s\le a_v\le b$.

Every member of $H$ lies below $a_t$. If $c\in C$ is an upper bound of $H$ with $a_t\nleq c$, apply the preceding paragraph to $b=a_t\wedge\neg c$. It supplies $h\in H$ with $0<h\le\neg c$, contrary to $h\le c$. Thus $a_t$ is the least upper bound in $C$.

In $D$ the join of $H$ is $j_C(t,n)$. If its nonzero remainder $r=a_t\wedge\neg j_C(t,n)$ belonged to $C$, the first paragraph applied to $b=r$ would give $0<h\le r$ for some $h\in H$. Yet $h\le j_C(t,n)$, a contradiction.
\end{proof}

The internal equality $\sup_C H=a_t$ holds at every candidate algebra. The next recursion seeks an algebra where the inclusion into $D$ preserves these particular joins. Each failure supplies its own new generator, the ambient remainder.

\begin{lemma}[Finite-term well-ordering]\label{lem:terms}
Let $\Gamma$ be a well-orderable set of generators for a Boolean algebra $C=\BA(\Gamma)$, where $\BA(\Gamma)$ denotes closure under finite Boolean operations. Then $C$ is well-orderable. Moreover, the set of finite histories over $C$ is well-orderable.
\end{lemma}
\begin{proof}
Fix a well-order of $\Gamma$, take a formal tagged copy of $\Gamma$, and adjoin a disjoint finite alphabet of Boolean operation symbols and punctuation. The resulting alphabet is well-ordered. Encode formal Boolean terms as finite words over this alphabet, ordered first by length and then lexicographically. For each fixed length, the lexicographic order on a finite power of a well-order is a well-order; the length-ordered sum over $\omega$ is therefore a well-order. Every element of $C$ is represented by at least one finite term. Assign to it the least term representing it. Distinct elements have distinct least terms, so this assignment injects $C$ into a well-ordered set and induces a well-order on $C$. Finally, finite sequences over a well-ordered set are well-orderable by length and then lexicographically. The histories form a subset of these finite sequences, so they are well-orderable as well.
\end{proof}

\begin{lemma}[Hartogs refinement]\label{lem:hartogs}
There is a well-orderable $C\in\mathcal C$ such that, for every $t\in L_C$ and $n<\omega$,
\begin{equation}\label{eq:ambient-sup}
a_t=\bigvee_D\{a_s:s\in J_C(t,n)\}.
\end{equation}
\end{lemma}
\begin{proof}
Let $\kappa=h(D)$ be the Hartogs ordinal of the underlying set of $D$, that is, the least ordinal which does not inject into $D$. We define a transfinite recursion up to $\kappa$, with a stop symbol propagated after the first successful stage.

Suppose no earlier stage has stopped and the elements $r_\beta$ have been defined for $\beta<\alpha<\kappa$. Let
\[
C_\alpha=\BA\bigl(\{d_{n,m}:n,m<\omega\}\cup\{r_\beta:\beta<\alpha\}\bigr).
\]
Label the displayed generators by the ordered sum of a fixed well-order of $\omega\times\omega$ and the ordinal $\alpha$. Different indices may name the same generator. Assign each generator its least index in this ordered sum; the resulting order is determined by the earlier sequence. By Lemma~\ref{lem:terms}, $C_\alpha$ and $L_{C_\alpha}$ are well-orderable. The construction in that lemma gives specified well-orders: in particular, the set $L_{C_\alpha}\times\omega$ has a fixed well-order determined by the earlier sequence.

If \eqref{eq:ambient-sup} holds for every pair $(t,n)\in L_{C_\alpha}\times\omega$, stop and record $C_\alpha$. Otherwise take the least failing pair $(t_\alpha,n_\alpha)$ and put
\[
r_\alpha=a_{t_\alpha}\wedge\neg j_{C_\alpha}(t_\alpha,n_\alpha).
\]
Since $j_{C_\alpha}(t_\alpha,n_\alpha)\leq a_{t_\alpha}$ and the pair fails the equality, $r_\alpha>0$. By Lemma~\ref{lem:residual}, $r_\alpha\notin C_\alpha$. In particular, $r_\alpha$ differs from every $r_\beta$ with $\beta<\alpha$.

This is a functional recursion: at each stage the least failing pair and the residual are uniquely specified. At a limit stage the generators are obtained from the set of all earlier $r_\beta$; Replacement ensures that this sequence and the generated subalgebra are sets. The order at a limit stage is computed anew from this sequence. No compatibility of the orders at different stages is needed. Thus the recursion is available in ZF. If it did not stop at any stage below $\kappa$, the map $\alpha\mapsto r_\alpha$ would be an injection from $\kappa$ into $D$, contrary to the definition of $h(D)$. Hence it stops at some $\alpha_*<\kappa$. The recorded algebra $C=C_{\alpha_*}$ is well-orderable and satisfies \eqref{eq:ambient-sup} for every history and every $n$.
\end{proof}

\begin{remark}
No monotonicity of the requirements is used. When the subalgebra grows, its history space and its cylinder sets change, so a requirement satisfied at one stage need not remain satisfied later. The recursion simply rechecks all requirements at every stage. Its termination depends only on the fact that each failure produces an element outside the algebra generated so far.
\end{remark}

\section{The complementing trees}

Fix a successful well-orderable algebra $C$ given by Lemma~\ref{lem:hartogs}. By Lemma~\ref{lem:terms}, $L_C$ is well-orderable. Fix an ordinal $\theta$ and a bijection
\[
\ell:L_C\longrightarrow\theta.
\]
We write histories in tree labels in place of their codes under $\ell$; formally, each occurrence of a history coordinate below means its image under $\ell$.

For $\varepsilon\in\{+,-\}$, define $T^\varepsilon\subseteq(\omega\times\theta)^{<\omega}$ as follows. A node of length $k$ is represented by
\[
\bigl((u(0),\ell(s_0)),\ldots,(u(k-1),\ell(s_{k-1}))\bigr),
\]
where $u\in\omega^k$ and $s_0,\ldots,s_{k-1}\in L_C$, subject to, for every $j<k$,
\begin{equation}\label{eq:tree-conditions}
s_{j-1}\pre s_j,\qquad
[s_j]_C\subseteq U_\varepsilon\cap O_j,\qquad
u\restriction(j+1)=u_{s_j}\restriction(j+1),
\end{equation}
where $s_{-1}=\varnothing$. Since each step is a proper extension, $|s_j|\geq j+1$, so the restriction of $u_{s_j}$ in \eqref{eq:tree-conditions} is defined. These conditions are inherited by initial segments, so $T^\varepsilon$ is a tree. Since all the data in \eqref{eq:tree-conditions} belong to the ground model, both trees belong to $V$. In particular the relation $[s]_C\subseteq U_\varepsilon\cap O_j$ is tested in $V$ when the node is admitted. A new branch in an extension only has to pass these recorded finite tests.

\begin{lemma}[Ground-model projections]\label{lem:ground-proj}
In $V$,
\[
p[T^\varepsilon]=A_\varepsilon\qquad(\varepsilon\in\{+,-\}).
\]
\end{lemma}
\begin{proof}
First suppose $(x,(s_j)_{j<\omega})$ is a branch through $T^\varepsilon$ in $V$. The histories are nested and $|s_j|\geq j+1$, so
\[
y=\bigcup_{j<\omega}s_j
\]
is an element of $X_C$. The coordinate requirements in \eqref{eq:tree-conditions} give $e_C(y)=x$. Moreover, $y\in[s_j]_C\subseteq U_\varepsilon\cap O_j$ for every $j$. By \eqref{eq:color-correct}, $x\in A_\varepsilon$. This proves $p[T^\varepsilon]\subseteq A_\varepsilon$.

For the reverse inclusion, fix $x\in R$ and begin with $t=\varnothing$, for which $\pi_x(a_t)=1$. Suppose a history $t$ has been chosen with $\pi_x(a_t)=1$, and let $n$ be the next stage. Applying the complete homomorphism $\pi_x$ to \eqref{eq:ambient-sup} yields
\[
1=\bigvee_{s\in J_C(t,n)}\pi_x(a_s).
\]
Since the join is taken in $2$, some $s\in J_C(t,n)$ satisfies $\pi_x(a_s)=1$. Choose the least such $s$ in the fixed well-order of $L_C$ and continue. This defines by recursion on $\omega$ a sequence
\[
s_0\pre s_1\pre\cdots,
\qquad \pi_x(a_{s_n})=1.
\]
Each selected cylinder is nonempty and is contained in one of the two color sets. Since the cylinders are nested, their colors must agree: if two successive colors differed, the later nonempty cylinder would be contained in two disjoint open sets. Let $\varepsilon$ be this common color.

For $j<|s_n|$, descent along the history gives
\[
a_{s_n}\leq d_{j,u_{s_n}(j)}.
\]
Applying $\pi_x$ and using $\pi_x(a_{s_n})=1$ shows that $x(j)=u_{s_n}(j)$. Consequently the sequence $(s_n)$, together with the real coordinate $x$, satisfies \eqref{eq:tree-conditions} and is a branch through $T^\varepsilon$. By the inclusion already proved, $x\in A_\varepsilon$. In particular, if $x$ was prescribed in $A_\delta$, disjointness of $A_+$ and $A_-$ forces $\varepsilon=\delta$, so the branch just constructed witnesses $x\in p[T^\delta]$.
\end{proof}

\begin{lemma}[Generic extension coverage]\label{lem:generic-cover}
For every $\mathbb B^+$-generic filter $G$ over $V$,
\[
V[G]\models\omega^\omega=p[T^+]\cup p[T^-].
\]
\end{lemma}
\begin{proof}
Work in $V[G]$ and fix $x\in\omega^\omega$. By Lemma~\ref{lem:all-reals}, choose $i\in I^V$ such that $x=i_G$. In Lemma~\ref{lem:ground-proj}, a join equal to $1$ in $2$ had a term equal to $1$. Here the join lives in $\mathbb B$: we must show that if its value belongs to $G$, then some term belongs to $G$. This is where genericity enters. Every history selected during the recursion is a finite tuple of elements of $C$, hence belongs to $V$. Fix such a $t\in L_C$ and suppose
\[
b=\pi_i(a_t)\in G.
\]
For a fixed $n<\omega$, applying $\pi_i$ in the ground model to \eqref{eq:ambient-sup} gives
\[
b=\bigvee_{s\in J_C(t,n)}c_s,
\qquad c_s=\pi_i(a_s).
\]
This equality uses only completeness of $\pi_i$. To pass from $b\in G$ to membership of some $c_s$ in $G$, define the ground-model subset of $\mathbb B^+$
\[
E_{i,t,n}=\{v\in\mathbb B^+:v\leq\neg b\text{ or }(\exists s\in J_C(t,n))\ v\leq c_s\}.
\]
We claim that it is dense. Given $q>0$, if $q\wedge\neg b>0$, take $v=q\wedge\neg b$. Otherwise $q\leq b$, and hence
\[
q=q\wedge b=q\wedge\bigvee_s c_s=\bigvee_s(q\wedge c_s).
\]
As $q>0$, some $q\wedge c_s$ is nonzero; take $v=q\wedge c_s$. In either case $v\in E_{i,t,n}$ and $v\leq q$.

The generic filter meets $E_{i,t,n}$. If $v\in G\cap E_{i,t,n}$, the alternative $v\leq\neg b$ is impossible because $b\in G$ and $G$ is a proper filter. Thus $v\leq c_s$ for some $s\in J_C(t,n)$, and upward closure of $G$ gives $c_s\in G$. We have proved
\begin{equation}\label{eq:generic-extension}
\pi_i(a_t)\in G\quad\Longrightarrow\quad
\exists s\in J_C(t,n)\ \bigl(\pi_i(a_s)\in G\bigr).
\end{equation}

Starting with $t=\varnothing$, for which $\pi_i(a_t)=1\in G$, apply \eqref{eq:generic-extension} successively for $n=0,1,\ldots$. At each stage choose the least available $s$ in the ground-model well-order of $L_C$. The parameters $i,t,n$ thus determine $J_C(t,n)$ and $E_{i,t,n}$ in $V$ at every finite stage, even when the resulting sequence of histories is new. We obtain
\[
s_0\pre s_1\pre\cdots,
\qquad \pi_i(a_{s_n})\in G.
\]
For each successive pair, the two finite histories and their recorded colors belong to $V$. There the later cylinder is nonempty and lies in the earlier cylinder, so their colors agree. Consequently every chosen history has the first color.

For $j<|s_n|$, we have $a_{s_n}\leq d_{j,u_{s_n}(j)}$, and therefore
\[
\pi_i(a_{s_n})\leq\pi_i(d_{j,u_{s_n}(j)})=i_{j,u_{s_n}(j)}.
\]
Since the left side belongs to $G$, so does the right side. By the definition of $i_G$, this implies
\[
u_{s_n}(j)=i_G(j)=x(j).
\]
Thus the histories form a branch through $T^\varepsilon$ over $x$. Since $x$ was arbitrary, the two projections cover all reals in $V[G]$.
\end{proof}

\begin{lemma}[Absolute disjointness]\label{lem:disjoint}
In every set-forcing extension of $V$, the projections of $T^+$ and $T^-$ are disjoint.
\end{lemma}
\begin{proof}
In $V$, form the synchronous product tree $S$ whose nodes are pairs $(v^+,v^-)$ of nodes of equal length from $T^+$ and $T^-$ with identical real coordinates. An infinite branch through $S$ in $V$ would give a real in $p[T^+]^V\cap p[T^-]^V$, contrary to Lemma~\ref{lem:ground-proj}. Thus $S$ has no infinite branch in $V$.

The node set of $S$ is well-orderable: it is a subset of the set of finite sequences over the well-orderable alphabet $\omega\times\theta\times\theta$. The proper-extension relation on $S$ is therefore well-founded. Indeed, if some nonempty subset $Y\subseteq S$ had a proper extension in $Y$ for each of its nodes, start with the least element of $Y$ in the fixed coding order and repeatedly take the least proper extension that remains in $Y$. This gives an infinite chain of strictly increasing lengths, whose union is a branch through $S$ in $V$, a contradiction.

By well-founded recursion there is an ordinal-valued rank $\operatorname{rk}_S$ on $S$ that strictly decreases along proper extensions. Explicitly $\operatorname{rk}_S(v)=\sup\{\operatorname{rk}_S(w)+1:v\pre w,\ w\in S\}$. The extension has the same finite nodes of $S$ and the same values of this ground-model function on them. Thus the strict rank inequalities persist in every forcing extension, where an infinite branch through $S$ would yield an infinite strictly decreasing sequence of ordinals. No such sequence exists. Hence $S$ has no branch in any set-forcing extension. If a real belonged to both tree projections in an extension, the two witnessing branches could be synchronized to form a branch through $S$. This is impossible.
\end{proof}

\begin{proof}[Proof of Theorem~\ref{thm:main}]
Assume (i), and fix a set forcing $\mathbb P$. The construction above gives ground-model trees $T^+,T^-$ on a common ordinal alphabet. Lemma~\ref{lem:ground-proj} gives their exact ground-model projections. Lemma~\ref{lem:generic-cover} gives forced coverage, and Lemma~\ref{lem:disjoint} gives forced disjointness. Thus they are $\mathbb P$-absolutely complementing, proving (ii).

Conversely, assume (ii). Then $A$ is universally Baire in the tree sense of Larson, Sargsyan, and Wilson. Their Lemma 2.4, proved in ZF, gives (i). This completes the proof.
\end{proof}

\section{Consequences and choice-free features}

\begin{corollary}[A single test space for each forcing]\label{cor:single-test}
For every set forcing $\mathbb P$ there are a space $X_{\mathbb P}$ with a clopen base and a continuous map $e_{\mathbb P}:X_{\mathbb P}\to\R$, both depending only on $\mathbb P$, such that for every $A\subseteq\R$, if $e_{\mathbb P}^{-1}[A]$ has the Baire property in $X_{\mathbb P}$, then there are $\mathbb P$-absolutely complementing ordinal trees whose ground-model projections are $A$ and $\R\setminus A$.
\end{corollary}
\begin{proof}
The space and map are those constructed in Section~2. The rest of the proof uses the Baire property of their single inverse image only to obtain the decomposition in Section~3. All subsequent constructions are carried out from that decomposition.
\end{proof}

\begin{corollary}[Complete ultrametric tests in ZF]
For $A\subseteq\R$, condition $(\ast)$ is equivalent to the Baire property of $f^{-1}[A]$ for every continuous $f:Y\to\R$ from a space admitting a compatible complete ultrametric. Completeness here means convergence of every Cauchy sequence. There is no restriction on density.
\end{corollary}
\begin{proof}
Open ultrametric balls are clopen, so $(\ast)$ implies the stated test. Conversely, fix $\mathbb P$. On $X_{\mathbb P}$ give distinct summands distance $1$. For distinct points $y,z\in X_C$, put $d(y,z)=2^{-k-1}$, where $k$ is their first coordinate of disagreement, and put $d(y,y)=0$. Agreement on initial segments gives the ultrametric inequality, and the tagged cylinders give the induced topology. A Cauchy sequence eventually stays in one summand $X_C$. For each $k<\omega$, sufficiently late terms have pairwise distance less than $2^{-k-1}$, so they agree on coordinates $0,\ldots,k$. The unique eventual values define a point $y\in C^\omega$. Each finite prefix of $y$ is a prefix of a sufficiently late term and hence belongs to $L_C$. Thus $y\in X_C$, and agreement on each finite prefix gives convergence to $y$. The test now applies to $e_{\mathbb P}$. Corollary~\ref{cor:single-test} gives the required trees for each $\mathbb P$, and Theorem~\ref{thm:main} gives $(\ast)$.
\end{proof}

\begin{remark}[Where choice is avoided]
The implication (i)$\Rightarrow$(ii) uses no form of the Axiom of Choice or Dependent Choice. After fixing one Baire decomposition, the selections in the recursions are canonical relative to fixed well-orders or to the natural order on $\omega$: the least natural number in Lemma~\ref{lem:extend}; the least finite term representing an element of a generated algebra; the least failing pair in the Hartogs recursion; and the least available history in the tree construction. The only global bound is the Hartogs ordinal of the set $D$. At no point is $D$, the set $I$ of all Boolean-valued real codes, or the collection of all forcing names well-ordered.
\end{remark}

\begin{remark}[No countable-intersection principle]
The proof does not invoke the Baire category theorem or assert that a countable intersection of dense open sets is nonempty. The topological input in the residual lemma is only that a finite intersection of dense open sets is dense. The countably many requirements are met one at a time through the ambient-supremum identity and, in a generic extension, through genericity for a ground-model dense set at each finite stage.
\end{remark}

\section*{Acknowledgements}
The author used ChatGPT and Codex to develop and assess candidate arguments for the Hartogs construction and the tree representation, to locate references, and to draft and edit the exposition. The author is responsible for the mathematical statements, proofs, and references in this paper.

\end{document}